\documentstyle[12pt]{article}
\newcommand{\be}{\begin{equation}}
\newcommand{\ee}{\end{equation}}
\newcommand{\ba}{\begin{eqnarray}}
\newcommand{\ea}{\end{eqnarray}}
\newcommand{\ban}{\begin{eqnarray*}}
\newcommand{\ean}{\end{eqnarray*}}
 
\renewcommand{\o}[2]{\frac{#1}{#2}}
\newcommand{\hf}{\o{1}{2}}
 
\newcommand{\qed}{\hspace*{\fill}\rule{3mm}{3mm}\quad}
\newcommand{\Pf}{\noindent {\em Proof.} }

\newcommand{\Rk}{\noindent {\em Remark.} }

\newcommand{\ra}{\rightarrow}
\newcommand{\pt}{\partial}

\newcommand{\Tr}{{\rm Tr}}
\newcommand{\tr}{{\rm tr}}

\newcommand{\te}{\tilde{e}}
\newcommand{\dM}{\partial M}

\newcommand{\End}{{\rm End}}

\newcommand{\Trs}{{\rm Tr}_s}

\newcommand{\Id}{{\rm Id}}
\newcommand{\hc}{\hat{c}}
\newcommand{\he}{\hat{e}}
\newcommand{\ext}{{\rm ext}}
\newcommand{\inn}{{\rm int}}
\newcommand{\Lap}{\Delta}
\newcommand{\cLap}{\Delta^c}

\font\BBb=msbm10 at 12pt 
\newcommand{\Bbb}[1]{\mbox{\BBb #1}}

\newtheorem{theo}{Theorem}[section]
 
\begin{document}
\newtheorem{lem}[theo]{Lemma}
\newtheorem{prop}[theo]{Proposition}  
\newtheorem{coro}[theo]{Corollary}
\title{\sc Analytic Torsion and R-Torsion for Manifolds with 
Boundary} 
\author{Xianzhe Dai\thanks{Department of Mathematics, University of California, Santa Barbara, CA 93106, U.S.A. (dai@math.ucsb.edu). Supported in part by NSF Grant DMS-9704296} \and 
Hao Fang\thanks{Department of Mathematics, Princeton University, Princeton, NJ 08540, U.S.A. (haofang@math.princeton.edu).}}
\date{}
\maketitle

\begin{abstract}  We prove a formula relating the analytic torsion and 
Reidemeister torsion on manifolds with boundary in the general case when the 
metric is not necessarily a product near the boundary. The product case has 
been established by W. Lu\"ck and S. M. Vishik. We find that the extra term 
that comes in here in the nonproduct case is the transgression of the Euler 
class 
in the even dimensional case and a slightly more mysterious term involving 
the second fundamental form of the boundary and the curvature tensor of the 
manifold in the odd dimensional case.
\end{abstract}

\section{Introduction}

The Reidemeister torsion (R-torsion for short) is a combinatorial/topological
invariant associated to a unitary representation of the fundamental group of a 
manifold. Introduced by
Reidemeister~\cite{r} and Franz~\cite{f}, it is used to classify the lens 
spaces. It was 
further developed by Milnor and Whitehead and used successfully in classifying 
the cobordisms. In searching for an analytic interpretation of the R-torsion, 
Ray 
and Singer \cite{rs1}, \cite{rs2} introduced the analytic torsion, which 
subsequently found many 
significant and far-reaching applications (Cf., for example, \cite{b}, \cite{bgs}, 
\cite{bl}, \cite{gs}, \cite{do}, \cite{w}). 

The celebrated Cheeger-M\"uller theorem \cite{c}, \cite{mu1} establishes the 
Ray-Singer conjecture: namely, 
on closed manifolds, the analytic torsion is indeed the same as the R-torsion.
Further significant work along this line includes that of M\"uller \cite{mu2} 
where he 
extends the theorem to unimodular representations, that of Bismut-Zhang 
\cite{bz}, 
which treats general representations (in which interesting secondary invariants come 
in), and that of Burghelea-Friedlander-Kappeler-McDonald \cite{bfkm} which deals with infinite dimensional representations (finite type Hilbert module). 

All these work deals with closed manifolds. With appropriate boundary 
conditions,
the torsion invariants can also be defined for manifolds with boundary. In 
fact, the torsions for manifolds with boundary provides an important stepping
stone in Cheeger's approach to the Cheeger-M\"uller theorem (Cf. \cite{c}, \cite{rs1}). In 
\cite{c}, the analog of Ray-Singer conjecture on manifolds with boundary is briefly discussed but the 
geometric information involved is not made explicit. This problem is addressed 
much later, first in Lott-Rothenberg \cite{lr}, and then L\"uck \cite{l} and 
Vishik \cite{v}. By 
assuming a product 
metric structure near the  boundary, they found that the analytic torsion and 
the 
R-torsion differ by a topological invariant, namely, (up to a constant), the 
Euler number of the boundary. This work has been generalized by Burghelea-Friedlander-Kappeler \cite{bfk} to the finite type Hilbert module.

We obtain a general formula of the difference of analyitc torsion and 
Reidemaster torsion on a Riemannian manifold with boundary, without assuming the product metric structure near the boundary. In this case, there are two extra, 
geometric terms coming in, involving the second fundamental form of the 
boundary.
In fact, one of the terms is the Chern-Simons class defined by Bismut-Zhang 
\cite{bz}, 
transgressing the Euler class. This term therefore vanishes if the dimension of 
the manifold is odd. The other term, which vanishes for the even dimension, is 
only slightly more complicated, involving 
the second fundamental form of the boundary and the curvature tensor of the 
manifold.

To be more precise, let $M$ be a compact Riemannian manifold with boundary $\dM$ 
and metric $g=g^{TM}$. If we denote $x$ the geodesic distance to the 
boundary $\dM$, then, near the boundary,
\[ g= dx^2 + g_{\dM}(x), \]
where $g_{\dM}(x)$ is a family of Riemannian metrics on the boundary. Let $g_0$ 
be a Riemannian metric which is of the form
\[ g_0=dx^2 + g_{\dM}(0) \]
near the boundary.
Let $\te(g_0, g)$ be the Chern-Simons class defined in \cite{bz}, which 
satisfies the equation
\be \label{transg}
d \te(g_0, g) = e(g) - e(g_0), 
\ee
where $e(g)$ denotes the Pfaffian form of $g$. Set, for $l\in [0, \ 1]$,
\be 
g_l=lg+(1-l)g_0.
\ee
Let $\rho:\, \pi_1(M) \longrightarrow U(N)$ be a unitary representation. Denote 
by $i$  the inclusion of $\dM$ into $M$. 

\begin{theo} \label{mt}
Let $T=T(M, \rho)$ be the analytic torsion of $(M, \ \rho)$ with respect to the
absolute 
boundary condition and $\tau=\tau(M, \rho)$ be the corresponding Reidemeister 
torsion. Then 
\[ \ln(T/\tau)= \chi (\dM, \rho)\ln 2 + {\rm rank}\,\rho\int_{\dM} i^* 
\te(g_0, g) + c\,{\rm rank}\,\rho \int_{\dM} i^*\phi, \]
where 
\[ c= \int_0^{\infty} \int_0^1 \int_0^{\infty} x(e^{-(x-x')^2/4s} + 
e^{-(x+x')^2/4s})
2 e^{-(x+x')^2/4(1-s)}\, dx'\, ds \, dx \]
and $\phi$ is a differential form on $M$ defined by the Berezin integral  
\[ \int^B \int_0^1  h_{ab} e^a
\wedge \hat{e}^b {\cal R}_0' e^{-{\cal R}}\, dl, \]
where $h_{ab}, \ 1 \leq a, b \leq n-1$, is the second fundamental form of the boundary, 
\[ {\cal R}_0'= \o{1}{4}R_{0jkl}\, e^j \wedge \hat{e}^k \wedge 
\hat{e}^l \]
and 
\[ {\cal R}=\o{1}{8} R_{ijkl}\, e^i 
\wedge e^j \wedge \hat{e}^k \wedge \hat{e}^l.\]
\end{theo}

\Rk One also obtain a similar formula for the relative boundary condition, with only 
sign changes.

\Rk The reader should consult \cite{ls} for examples where the these geometric contributions
are shown to be nonzero, using the formulas given in \cite{bg}. The article \cite{bg} also contains a
thorough discussion on the aboslute and relative boundary conditions.
 
An immediate application of this result is a gluing formula for the analytic 
torsion, since the gluing formula for R-torsion is well known. 

Our strategy of proof is a very natural one. Namely we deform to the product case, where the difference 
of the analytic torsion and R-torsion is known by \cite{l}, \cite{v}. The variation of 
the log of the ratio of the analytic torsion over the R-torsion is given by the 
constant term in certain asymptotic expansion of the heat kernel. To compute 
this term we employ the local index theorem technique, much in the spirit of
Bismut-Zhang \cite{bz}. The new difficulty, however, is that the boundary condition introduces certain non-uniform 
behavior in the local heat asymptotic expansion, in the form of some Gaussian type concentration near the boundary. This difficulty is overcome by 
bringing in the technique of Melrose (Cf. \cite{dm}). The non-uniform behavior 
of the heat kernel for the absolute boundary condition (say) is resolved by 
lifting it to a larger space obtained by performing certain blowup operation on the usual carrier 
of the heat kernel.  This technique also effectively separates the interior
and boundary contributions. The interior contribution is zero, just as in the closed case,
while the boundary contribution comes from a model problem on the half tangent space of the
manifold at the boundary.

\medskip
{\sc Acknowledgment}.  The authors are very grateful to Xiaonan Ma for his careful reading of an earlier version and for pointing out several mistakes and misprints. The authors would like to thank Jean-Michel Bismut, Alice Chang, Paul Yang and Weiping Zhang for their interests and useful discussion. The second author would also like to thank Alice Chang, Paul Yang and Weiping Zhang for their constant encouragement.

\section{Analytic torsion for manifold with boundary}

In this section we recall the definition of the absolute and relative boundary conditions,
their corresponding torsion invariants, while fixing our notations. We also state the variational 
formula for the log of the ratio of the analytic torsion over the R-torsion (from \cite{rs1} and \cite{c}), which is our starting point.

\subsection{Absolute and relative boundary conditions}

Let $M$ be a compact  manifold with boundary and $g$ be 
a Riemannian metric on $M$. If we denote $x$ the geodesic distance to the 
boundary $\dM$, then, near the boundary,
\[ g= dx^2 + g_{\dM}(x), \]
where $g_{\dM}(x)$ is a family of Riemannian metrics on the boundary. Let $\xi 
\ra M$ be the flat bundle associated to a representation $\rho:\, \pi_1(M) \ra 
U(k)$. Consider the Hodge Laplacian 
\be \label{hlap}
\Delta=(d+\delta)^2:\ C^{\infty}(M; \Lambda(M)\otimes \xi) \ra C^{\infty}(M; 
\Lambda(M)\otimes \xi). 
\ee

At the boundary we have the splitting
\be \label{spff}
\Lambda(M)\otimes \xi |_{\pt M}=\Lambda(\pt M)\otimes \xi\oplus\Lambda(\pt 
M)\otimes \xi
\ee
corresponding to the decomposition for a form $\omega \in C^{\infty}(M; 
\Lambda(M)\otimes \xi)$:
\[ \omega=\omega_t+dx\wedge\omega_n, \  \omega_t,\, \omega_n \in C^{\infty}(M; 
\Lambda(\dM)\otimes \xi) \]
near the boundary. Define a linear map $\sigma$:
\[ \sigma(\omega)=\omega_t-dx\wedge\omega_n. \]
Then $\sigma$ is self adjoint and $\sigma^2=1$. Moreover the splitting 
(\ref{spff}) corresponds to the decomposition into the $\pm 1$-eigenspace of 
$\sigma$. 

{}From the splitting we define two projections
\ban
 P_a: \  C^{\infty}(\pt M; \Lambda(M)\otimes \xi |_{\pt M}) & \ra & 
C^{\infty}(\pt M; \Lambda(\pt M)\otimes \xi),  \\
 P_a(\omega) & = & \omega_n|_{\pt M}; \\
P_r:\  C^{\infty}(\pt M; \Lambda(M)\otimes \xi |_{\pt M}) & \ra & 
C^{\infty}(\pt M; \Lambda(\pt M)\otimes \xi),  \\
 P_r(\omega) & = & \omega_t|_{\pt M}.
\ean
i.e., $P_a$ is the orthogonal projection onto the $-1$-eigenspace of $\sigma$ 
and $P_r$ the orthogonal projection onto the $+1$-eigenspace. The absolute and 
relative boundary conditions for the Hodge Laplacian are
\ban \label{bc}
{\rm absolute}: &  \left\{ \begin{array}{l} P_a(\omega|_{\pt M})=0, \\
 P_a(d\omega|_{\pt M})=0. \end{array} \right. \\ 
{\rm relative}: &  \left\{ \begin{array}{l} P_r(\omega|_{\pt M})=0, \\
 P_r(\delta\omega|_{\pt M})=0. \end{array} \right.
\ean

It is well known that these define elliptic boundary conditions for the Hodge 
Laplacian, with the corresponding cohomology being the absolute (resp. relative) cohomology of $M$. Hence the names. 

To get a flavor of these boundary conditions, the special case when $M=\dM \times {\Bbb R}_+$, and $g=dx^2+g_{\dM}$, 
where $x$ is the variable of ${\Bbb R}_+$ and $g_{\dM}$ is a metric on the 
boundary, deserves some elaboration. Using separation of variables, we can write 
\be 
\omega = f(x) \omega_1^{\dM} + g(x) dx\wedge \omega_2^{\dM},
\ee
where $ \omega_1^{\dM}$ and $ \omega_2^{\dM}$ are differential forms on $\dM$.
Then the absolute and relative boundary conditions are 
reduced to a combination of Dirichlet and Neumann boundary conditions:
\ba \label{bc}
{\rm absolute}: &  \left\{ \begin{array}{l} g(0)=0, \\
\o{\pt f}{\pt x}(0)=0. \end{array} \right. \\ 
{\rm relative}: &  \left\{ \begin{array}{l} f(0)=0, \\
 \o{\pt g}{\pt x}(0)=0. \end{array} \right.
\ea
This greatly simplifies our computation of the model problem.

Back to the general discussion, it follows easily from the Stokes' Theorem that
\be \qquad
(\Delta\omega, \theta )=(\omega,\Delta \theta)\! +\!\int_{\pt M}\omega\wedge
*d\theta\! -\!\int_{\pt M}\theta\wedge *d 
\omega \!+\!\int_{\pt M}d^*\omega\wedge
*\theta \!-\!\int_{\pt M}d^*\theta\wedge * \omega .
\ee
The boundary terms all vanish if both $\omega$, $\theta$ satisfy the absolute
(resp. relative) boundary condition. i.e. the absolute and  relative boundary conditions
 are also self adjoint boundary 
conditions. 

Therefore, if $F_1(t, x, z)$, $F_2(t, z, y)$ are double forms on $M$ that 
both satisfy 
the absolute (resp. relative) boundary condition, then (Cf. \cite[p275]{c})
\begin{eqnarray*}
\lefteqn{(F_1(0,x, \cdot),\, F_2(t, \cdot, y))- (F_1(t,x, \cdot),\, F_2(0, 
\cdot, y))} \nonumber \\
& & =\int_0^t (F_1(t-s,x, \cdot),\, (\pt_t +\Delta)F_2(s, \cdot, y)) ds
- \int_0^t ((\pt_t +\Delta) F_1(t-s,x, \cdot),\, F_2(s, \cdot, y)) ds.
\end{eqnarray*}
Thus if we denote 
\be F_1*F_2 = \int_0^t \int_M \langle F_1(t-s,x, z),\, F_2(s, z, y)\rangle dz\, 
ds,
\ee
then we have arrived at the so called Duhamel principle:

\begin{prop} Let $e^{-t\Delta}$ denote the heat kernel of the Laplacian with the 
absolute (resp. relative) boundary condition. Let $F_1(t, x, y)$ be a double 
form 
on $M$ which satisfies the absolute (resp. relative) boundary condition, and is 
the delta function at $t=0$. Then
\be \label{dp}
e^{-t\Delta}= \sum_{m=1}^{\infty} (-1)^{m-1} F_m ,
\ee
where $F_m= F_{m-1}*(\pt_t +\Delta) F_1$.
\end{prop}

\Pf This follows by iterating the formula above.
\qed

\subsection{Analytic torsion and R-torsion}

With elliptic self adjoint boundary conditions at our disposal the analytic torsions
for manifolds with boundary can be defined using exactly the same formula as in the closed case.
We only need to replace the usual heat kernel with the heat kernel corresponding to the appropriate 
boundary conditions. For simplicity it will be implicitly assumed that our heat kernels correspond to
the absolute boundary condition, unless stated otherwise.

Thus, the analytic torsion $T(M, \rho)$ is defined in terms of the torsion zeta function
\be
\log T(M,\rho)=\zeta'_T(0) \label{dat} 
\ee
where $\zeta_T(s)$ is a meromorphic function of $s\in {\Bbb C}$, defined as follows. For $\Re\, s>>0$
\be \label{tzeta}
\zeta_T(s)=\frac{1}{\Gamma(s)}\sum_{p=1}^{\dim M}p (-1)^{p+1}
\int_{0}^{\infty}t^s\Tr\left(e^{-t\Delta_p} P^{\perp}\right)\frac{dt}{t},
\ee
where $\Delta_p$ is the Laplacian on $p$-forms satisfying the absolute boundary 
conditions, and $P^{\perp}$ denote the projection onto the orthogonal complement 
of the harmonic forms.  Using the standard elliptic theory, one sees that $\zeta_T(s)$ extends to 
 a meromorphic function of $s\in {\Bbb C}$ which is regular 
at $s=0.$ In particular, the right hand side of (\ref{dat}) is well-defined.

The natural ${\Bbb Z}$ grading on the space of differential forms induces a natural  ${\Bbb Z}_2$
grading. We denote the corresponding supertrace $\Trs$. With this notation formula (\ref{tzeta})
can be rewritten as
\[ \zeta_T(s)=- \frac{1}{\Gamma(s)}\int_{0}^{\infty}t^s\Trs\left(N e^{-t\Delta} P^{\perp}\right)\frac{dt}{t},
\]
where $N$ denotes the so called number operator which mutiplies a differential $p$-form by $p$.

Similarly, the corresponding R-torsion $\tau(M, \rho)$ can be defined the same way as in the closed 
case using the absolute cochain complex.

The analytic torsion as defined above depends on the Riemannian metric $g$. So does
the R-torsion. The variation of the torsion invariants with respect to the metric change
is well studied. We quote the following result from \cite{rs1}, \cite{c}.
\begin{prop}
Let $g_l$ be  a family of Riemannian metrics ($l\in [0,1]$),  and
$*_l$  the corresponding Hodge $*$-operator. Denote $T_l=T_l(M, \rho)$ (resp. $\tau_l=\tau_l(M, \rho)$) the analytic torsion (resp. R-torsion) with respect to $g_l$. If the unit normal of $\dM$ with respect to
$g_l$ is independent of the parameter $l$, then
\be \label{vf}
\o{d}{dl} \log T_l/\tau_l= \lim_{t\ra 0}\Trs (*_l^{-1} 
\o{d*_l}{dl}e^{-t\Delta_l})
\ee
\end{prop}

\Rk In general the right hand side of (\ref{vf}) should be interpreted as the constant term in the 
asymptotic expansion of $\Trs (*_l^{-1} 
\o{d*_l}{dl}e^{-t\Delta_l})$. In our case, however, all the negative powers of 
the asymptotic expansion drop out, as we will see from our proof, and the limit
therefore is the genuine limit.

Thus, to prove our result, we deform to the product case. Namely, if we denote 
$x$ the geodesic distance to the boundary $\dM$ with respect to the Riemannain metric $g$ on 
$M$, then $g$ can be written as
\be 
g=dx^2 + g_{\dM}(x)
\ee
near the boundary. Now, let $g_0$ be a Riemannian metric which is of the form
\be
g_0=dx^2 + g_{\dM}(0)
\ee
near the boundary. Set 
\be 
g_l=lg+(1-l)g_0. \label{md}
\ee
We compute the right hand side of (\ref{vf}), using the local index theorem 
technique in the spirit of \cite{bz}. However, before that is possible, we first need to develop a thorough 
understanding of the uniform structure of the heat kernel for absolute (or relative) boundary condition.

\section{The Uniform Structure of Heat Kernel}

In this section we analyze the uniform structure of the heat kernel for 
elliptic boundary value problems of the type discussed above. The local index 
theorem technique dictates that we study the pointwise trace of the heat kernel. 
However, the elliptic boundary conditions introduce certain nonuniform behavior 
in the pointwise trace of the heat kernel in the form of some Gaussian type 
concentration at the boundary. We resolve this nonuniform behavior by studying 
the heat kernel in a bigger space, by introducing certain singular coordinates.
In fact, we will construct a pseudodifferential calculus adapted to the 
situation and show that the heat kernel is certain element in this calculus.
All the regularity and uniform property is contained in this statement.

\subsection{Construction of the pseudodifferential calculus}

Let $Z=M\times M \times {\Bbb R}_+$ be the usual carrier for the heat kernel. If 
we denote by $(x, y)$ the local 
coordinates of $M$ near $\dM$, then $(x, y, x', y', t)$ will be a local 
coordinates on $Z$. Our analysis of the heat kernel will be carried out in a 
blowup version of the space $Z$. For a space $X$, we denote by $\Delta(X)$ the 
diagonal in $X\times X$. Let $B_1=\Delta(\dM) \times \{0\}$
and $B_2=\Delta(M) \times \{0\}$ be subspaces of $Z$. Also denote by $S$ the 
parabolic bundle generated by $dt$. Consider the space \cite{dm}
\be \label{hds}
M_{h0}^2=[Z; B_1, B_2, S]
\ee
obtained from $Z$ by blowing up, parabolically in $t$, first $B_1$, then 
$B_2$. The blow down map is denoted by 
$\beta_{h0}$. That is, $\beta_{h0}:\, M_{h0}^2 \ra Z$.

The new space $M_{h0}^2$ is a manifold with corner. For the moment let's look at 
the structures of the 
boundary hypersurfaces. There are three of them lying above 
$\{t=0\}$. Among these we have 
$B_{f\! f}$ from the blow up of $B_1$, called the front face; 
$B_{t\! f}$ from the blow up of $B_2$, called the temporal face. These two are 
fibered over the 
submanifolds that are blown up. 
In fact $B_{t\! f}$ can be viewed as a natural compactification of the 
rescaled tangent bundle $^{0}TM$ of $M$, whose sections are precisely vector fields vanishing at the boundary, as can be verified using suitable coordinates such as (\ref{pc2}). Similarly 
$B_{f\! f}$ is related to the half tangent bundle of $M$ at $\dM$ as follows: the front face $ B'_{f\! f}$ resulting from the first blowup fibers over $B_1\cong \pt M$, and in fact is canonically isomorphic to
\be \label{htb}
{\Bbb R}_+ \times T_{\dM}^+M  \ra \pt M,
\ee
with $T_{\dM}^+M$ denoting the half tangent bundle of $M$ at $\dM$. Then $B_{f\! f}$ is a suitable blowup of $B'_{f\! f}$, namely, if $T$ is the coordinate of ${\Bbb R}_+ $, $(s, \bar{Y}, y')$ the fiber coordinates of $ T_{\dM}^+M$, then $B_{f\! f}$ is obtained by fiberwise blowing up $T=0, s=1, \bar{Y}=0$. 

The rest of the boundary hypersurfaces arises from the lift of those of 
$Z$ under the blowdown map. More precisely we have the lift of the hypersurface $t=0$ of $Z$, denoted by 
$B_{tb}(M_{h0}^2)$; the lift of the hypersurface $x=0$ of $Z$, $B_{lb}(M_{h0}^2)$; 
and 
the lift of $x'=0$, $B_{rb}(M_{h0}^2)$. 
Note that $B_{t\! f}$ only meets $B_{f\! f}$, $B_{tb}$.

In effect, the blowup is equivalent to introducing certain singular coordinates near the boundary. 
In the first blowup, the projective coordinates
\be \label{pc1}
T=t/x'^2, \  s=x/x', \  x', \bar{Y}=\o{y-y'}{x'}, \  y'
\ee
will be valid coordinates near the front face, except at the corner with 
$B_{rb}$. After the second blowup, we can 
use the projective coordinates
\be \label{pc2}
T, \  X=\o{s-1}{\sqrt{T}}=\o{x-x'}{\sqrt{t}}, \  
Y=\o{\bar{Y}}{\sqrt{T}}=\o{y-y'}{\sqrt{t}}, \ x', \  y'
\ee near the temporal face, except at the corner with $B_{tb}$.

We will denote by $\rho$ the defining function for each boundary hypersurface, 
with 
the appropriate subscript attached. Thus $\rho_{f\! f}$ will denote the defining 
function of $B_{f\! f}$, etc.

The pseudodifferential calculus will consists of operators whose kernels are 
defined on the space $M_{h0}^2$ and normalized with respect to 
the half-density
\be \label{hfdb}
 K\! D_{h0}=\rho_{f\! f}^{-(n+2)/2}\rho_{t\! f}^{-(n+3)/2}\Omega^{1/2}(M_{h0}^2), 
\ee
where $\Omega^{1/2}$ denotes the standard half-density bundle.
Let $I$ denote the index set $\{k_1, k_2, k_3, k_4\}$. The space of operators 
with index set $I$ is defined to be
\[ 
\Psi_{h0}^{I}(M, \Lambda M\otimes \xi)=\rho_{f\! f}^{k_1}\rho_{t\! 
f}^{k_2}\rho_{lb}^{k_3}\rho_{rb}^{k_4}\rho
_{tb}^{\infty}C^{\infty}(M_{h0}^2; {\rm Hom}(\Lambda M\otimes \xi) \otimes
K\! D_{h0}). \]

To see how the operators in $\Psi_{h0}^{I}$ act, consider the bilinear map
from
\[\dot{C}_0^{\infty}([0,\, \infty) \times M;\Lambda M\otimes \xi\otimes 
\Omega^{1/2})\times \dot{C}^{\infty}([0,\, \infty) \times M;\Lambda M\otimes 
\xi\otimes \Omega^{1/2})\]
 to 
\[ \dot{C}^{\infty}(Z;\Lambda M\otimes \xi\otimes \Omega^{1/2}) \]
 defined by
\be
 \phi \hat*_t \psi = \int\limits_0^\infty\phi(t+t',x)\psi(t',x')dt' 
\ee 
Lifting to $M_{h0}^2$ we can define
\be
\langle A\psi,\phi\rangle=
\int\limits_{M_{h0}^2}A\cdot\beta_{h0}^*(\phi\hat*_t\psi).
\ee

It can be shown that these operators form a filtered algebra containing the heat 
kernel of an elliptic boundary value problem. However we do not need the full 
strength of this result and therefore will only be establishing those that we need in 
the following subsections.

\subsection{Normal homomorphisms}

Just like the classical pseudodifferential operators, the regularity property of 
elements of $\Psi_{h0}^{I}$ is described by their symbols, except that there are 
more than just one (principal) symbol here, called the normal homomorphisms. 

The normal homomorphism at $B_{t\! f}$, or the heat homomorphism, is defined by 
dividing out $T^{\o{k_2}{2}}$ 
and restricting to $B_{t\! f}$:
\be 
N_{h, k_2}:\ \Psi_{h0}^{I} \ra \rho_{f\! f}^{k_1} \rho_{lb}^{k_3} 
\rho_{rb}^{k_4}\rho_{tb}^{\infty}C^{\infty}(B_{tf}; K\! D_{h0}|_{B_{t\! f}}). 
\ee
Since elements of $\Psi_{h0}^{I}$ vanish rapidly at $B_{tb}$ the heat homomorphism simply describes the leading term at $B_{t\! f}$. We note that $B_{t\! 
f}$ meets only $B_{f\! f}$ and $B_{tb}$. Moreover, since $B_{t\! f}$ is the 
natural compactification of $^{0}TM$, $K\! D_{h0}|_{B_{t\! f}}$ is 
canonically isomorphic to the fiber density 
bundle of $^{0}TM$ and thus  the heat homomorphism can be rewritten as
\[ N_{h, k_2}:\ \Psi_{h0}^{I} \ra 
\rho_{f\! f}^{k_1}{\cal S}(^{0}TM; 
\Omega_{\mbox{fiber}}). \]
Here ${\cal S}(^{0}TM; \Omega_{\mbox{fiber}})$ denotes the space of smooth fiber
densities which are of Schwarz class along fibers of $^{0}TM$.

The normal homomorphism at $B_{f\! f}$ describes the leading term at $B_{f\! 
f}$. Thus we divide out 
$\rho_{f\! f}^{k_1}$
and restrict to $B_{f\! f}$:
\be \label{nhff} N_{f, k_1}:\ \Psi_{h0}^{I} \ra 
\rho_{t\! 
f}^{k_2}\rho_{lb}^{k_3}\rho_{rb}^{k_4}\rho_{tb}^{\infty}C^{\infty}(B_{f
f}; K\! D_{h0}|_{B_{f\! f}}). 
\ee
By the previous discussion on the structure of $B_{f\! f}$ (see (\ref{htb}))  
the range of $N_{f, k_1}$ is isomorphic to the space of smooth functions on 
${\Bbb R}_+ \times T_{\dM}^+M$,  except for conormal concentration along the submanifold $s=1, \bar{Y}=0$ 
at $T=0$, with the order of vanishing at the boundaries 
$s=0$ and $s=\infty$ prescribed by the index set, and  vanishing rapidly at 
$\bar{Y}$-infinity.
In other words, the range of $N_{f, k_1}$ is the (restricted) fiberwise heat calculus on 
$T_{\dM}^+M$.

For an operator in $\Psi_{h0}^{I}$, its values under the normal homomorphisms 
will be called its normal operators. Individually, each normal homomorphism is 
surjective. However the normal 
operators for an element of $ \Psi_{h0}^{I}$ have to agree at the common 
corners. These are the compatibility conditions. On the other hand, since 
essentially just smooth functions are involved, the compatibility conditions 
are clearly the only obstructions to the existence of operators with the
prescribed normal operators.This fact will be used in constructing heat kernels 
as an element of $\Psi_{h0}^I$.

\subsection{Compositions}

Another ingredient in constructing heat kernels as an element of $\Psi_{h0}$ 
will be the composition properties of elements of $\Psi_{h0}$.
We first study the composition of these operators with differential 
operators. 

\begin{prop} \label{comp1}
Let $A \in \Psi_{h0}^I (M, \Lambda M\otimes \xi)$. 
If $V$ is any smooth vector field on $M$, and $T=t/x'^2$,  
then we have
\[ (t^{\hf} V) \circ A \in \Psi_{h0}^I(M, \Lambda M\otimes \xi) \]
and
\be \label{nop1}
\qquad N_{f,k_1} (t^{\hf}V \circ A) = T^{\hf} \sigma_1 (V) N_{f,k_1} (A), \ \ 
N_{h,k_2} (t^{\hf}V \circ A) = \sigma_1 (V) N_{h,k_1} (A). 
\ee
Here $\sigma_1 (V)$ denotes the usual symbol of a vector field. 
Furthermore, 
\[ (t\partial_t \circ A) \in  \Psi_{h0}^I(M, \Lambda M\otimes \xi) \]
and 
\be \label{nop2}
\qquad N_{f,k_1}(t\partial_t\circ A)=T\pt_T N_{f,k_1}(A), 
N_{h,k_2}(t\partial_t\circ A)=-\frac{1}{2}(R+n+2+k_2)N_{h,k_2}(A).
\ee
\end{prop}

\Pf We first assume that $\rho$ is the trivial representation, i.e. $\xi$ is  
the trivial bundle. Clearly we can restrict our attention to a region near the 
front and temporal faces. Therefore we assume that $A$ is supported near such a 
region. Further, by partition of unity, we can assume that $A$ is either 
supported near the front face but away from the temporal face, or is supported 
near the temporal face. 

\def\hsp{\hspace{-1pt}}
For the first case, we get to use the projective coordinates (\ref{pc1}),
But it is more convenient to use the project coordinate
\be \label{mcpc}
t, \ X=\frac{x}{t^{\hf}}, \ X'=\frac{x'}{t^{\hf}}, \ 
Y=\frac{y-y'}{t^{\hf}}, \ y',
\ee
which is valid near the front face, except at the corner,
where $X, X', Y=\infty$. In this coordinate 
$A\! \in\! \Psi_{h0}^I(\hsp M \hsp)$
can be written as $A\! =\! t^{\frac{k_1}{2}-\frac{n+2}{4}} 
a(t,X, X', Y, y')  |dt dX dX' dY dy'|^{\hf}$,
where $a$ is a smooth function of $(t^{\hf},X, X', Y, y')$ and 
vanishes rapidly as $X, X', Y \rightarrow \infty$. In this representation
$N_{f, k_1}(A)=a(0, X, X', Y, y') |dXdX'dYdy'|^{\hf}$.

Now if $\phi=\phi_0|dtdxdy|^{\hf}$ and $\psi=\psi_0 |dtdxdy|^{\hf}$ we have 
\[ \langle A\psi,\phi\rangle=
\int_{M_h^2} t^{\frac{k_1}{2}} a \beta_h^*(\phi_0\hat*_t\psi_0)
 |dt dX dX' dY dy'|.   \]

Also, if we let $V=v_1(x,y)\pt_x+ v_2(x, y)\pt_y $ be a smooth vector field and 
$V'$ denote its transpose:
\[ \langle V\psi, \phi\rangle= -\langle \psi, V'\phi\rangle, \]
then $V'=\pt_x\circ v_1+\pt_y \circ v_2$,  and $\beta_h^*(t^{\hf}V')=t^{\hf}\pt_x(v_1) + 
t^{\hf}\pt_y(v_2) + v_1\pt_X + v_2\pt_Y$.

A computation similar to that of \cite[p44]{dm} then
shows that
\[ (t^{\hf}V) \circ A = t^{\frac{k_1}{2}-\frac{n+2}{4}} 
(v_1 \pt_X(a) + v_2 \pt_Y(a) + t^{\hf}(\pt_xv_1) a+ t^{\hf}(\pt_yv_2) a)|dt dX dX' dY 
dy'|^{\hf} \]
and therefore $(t^{\hf}V) \circ A \in \Psi_{h0}^I(M)$, with
\be \label{nopipc}
 N_{f,k_1} (t^{\hf}V \circ A) = (v_1 \pt_X + v_2 \pt_Y) a(0, X, X', Y, y') 
|dXdX'dYdy'|^{\hf}.
\ee
Now we change the coordinate from (\ref{mcpc}) to (\ref{pc1}) and get
\be
\pt_X =T^{\hf}\pt_s, \hspace{1in} \pt_Y= T^{\hf}\pt_{\bar{Y}},
\ee
which, when plugged into (\ref{nopipc}), immediately gives 
\[ N_{f,k_1} (t^{\hf}V \circ A) =T^{\hf}\sigma_1 (V) N_{f,k_1} (A).\]

Similarly we have $t\pt_t \circ A \in  \Psi_{h0}^I (M)$ and
\be \label{noftd}
N_{f,k_1} (t\pt_t \circ A) = [\frac{k_1}{2}  -\frac{1}{2}(n+1+ X\pt_X+ X'\pt_X'+ 
Y\pt_Y)] N_{h,-j} (A).
\ee
But changing to the coordinate (\ref{pc1}) we see 
\[ N_{f,k_1} (t\pt_t \circ A) = T\pt_T N_{f,k_1} (A). \]

For the second case when the support of $A$ is near the temporal face, the 
normal operators can be obtained in the same way using the projective 
coordinates
(\ref{pc2}).

Finally, if $\rho$ is not trivial, by linearity, we can assume that $\phi,$ 
$\psi$ and $A$ are supported in a small neighborhood where $\xi$ is trivialized 
by an orthonormal basis $\{s_i\}$. Write 
\[ \phi=\phi_i s_i, \ \psi=\psi_i s_i, \ A=A_{ij}s_i^*\otimes s_j. \]
Then
\ban
\langle (V\circ A)\psi, \phi \rangle & = & -\langle A\psi, \nabla_V\phi \rangle
\\ & = & -\langle A_{ij}\psi_i, V\phi_j+\Gamma_{kj}(V)\phi_k\rangle,
\ean
where $\Gamma_{kj}(V)=\langle \nabla_Vs_i, s_j \rangle$. This reduces to the 
scalar case and the connection produces only a lower order term. 
\qed

We now consider the composition of elements of $\Psi_{h0}^{I}(M; \Lambda M\otimes \xi)$. For 
simplicity we will not consider the full composition properties, only the 
composition of the following Voltera type operators, which suffices for our 
purpose.
Thus if ${\cal E}=(E_{f\! f},E_{lb}, E_{rb})$ is an index family for
$M^2_{h0},$ assumed trivial at $B_{t\! f}$ and $B_{tb}$ then let
\[ \Psi_{h0}^{-\infty,{\cal E}}(M,\Omega^{\hf})= \Psi_{h0}^{E_{f\! f}, \infty, 
E_{lb}, E_{rb}}(M, \Lambda M\otimes \xi) \]
be the space of polyhomogeneous conormal distributions on $M_{h0}^2$ which 
vanish
rapidly at $B_{t\! f}$ and $B_{tb}$ and have expansions at $B_{f\! f}$, $B_{lb}$ 
and $B_{rb}$ with
exponents from $E_{f\! f}$, $E_{lb}$ and $E_{rb}$ respectively.

\begin{prop} \label{comp2}
We have
\[ \Psi_{h0}^{-\infty,{\cal E}}(M, \Omega^{\hf})\circ
\Psi_{h0}^{-\infty,{\cal F}}( M, \Omega^{\hf})\subset
\Psi_{h0}^{-\infty,{\cal G}}( M, \Omega^{\hf}) \]
where
\ban
G_{f\! f} & = & \left[E_{f\! f} + F_{f\! f} \right] \cup \left 
[E_{lb}+F_{rb}+n+2\right] \\
G_{lb} & = & \left[E_{f\! f} +F_{lb} \right] \cup \left[E_{lb}\right] \\
G_{rb} & = & \left[F_{f\! f} +E_{rb} \right] \cup \left[F_{rb}\right].
\ean
\end{prop}

\Pf  Since the operators here all have their Schwartz kernel vanishing rapidly at $B_{t\! f}$, we might as well consider them as living in
\[ M_0^2=[Z; B_1, S]. \]

 The proof relies on lifting everything to a triple space. For this reason we consider 
\[ W={\Bbb R}_+^2 \times M^3 =\{ (t, t', m, m', m'')\ | \ t\geq t'\}, \]
and the following submanifolds of $W$:
\[ B_F =\{ (t, 0, m, m', m') \}, \ \ B_S= =\{ (t, t, m, m, m'') \}, \ \
B_C =\{ (0, t', m, m', m) \}, \]
together with
\[ B_T =\{ (0, 0, m, m, m) \}, \]
which is the intersection of the all three. Let $S_F, S_S, S_C$ be the lift of the parabolic bundle $S=\{ dt\}$ to the corresponding submanifolds $B_F, B_S, B_C$ and $S_T=\{dt, dt'\}$ the parabolic bundle on $B_T$. Then the triple space $W_{0}$
is obtained from $W$ by appropriate blowups:
\[ W_{0}=[W; B_T, S_T; B_F, S_F; B_S, S_S; B_C, S_C]. \]

Once again we obtain a manifold with corner. We label the hypersurfaces created by the blowups by $f\! F, f\! f, s\! f, c\! f$ respectively, and those from the lift by $fb, sb, cb$. The triple space $W_0$ is related to the double space $M_0^2$ by the $b$-fibrations
\[ \pi_F,\ \pi_S, \ \pi_C: \ W_0 \ra M_0^2, \]
corresponding to the projections $p_F, \ p_S, \ p_C:\ W \ra Z$,
\ban 
 p_F(t, t', m, m', m'') & = & (t', m', m''), \\
 p_S(t, t', m, m', m'') & = & (t-t', m, m'), \\
 p_C(t, t', m, m', m'') & = & (t, t', m, m''). 
\ean

Let $f:\ X\longrightarrow Y$ be a $b$-map between manifolds with corners,
i.e.,  a $C^{\infty}$ map such that if $\rho_i'\in C^{\infty}(Y),$ $i=1,\dots,N'$ are defining
functions for the boundary hypersurfaces of $Y$ and $\rho_j\in C^{\infty}(X),$
$j=1,\dots,N$ are defining functions for the boundary hypersurfaces of $X$
then
\[ f^*\rho_i'=a_i\prod\limits_{j=1}\rho_j^{k(i,j)},\ 0<a_i\in C^{\infty}(X).\]
The non-negative integers $k(i,j)$ are the boundary exponents of $f.$ In the following table we computed all boundary exponents of the $b$-fibrations $\pi_F,\ \pi_S, \ \pi_C$.

\medskip
\begin{center}
\begin{tabular}{|l|l||l|l|l|l|l|l|l|l}      \hline
 & & fF & ff & sf & cf & fb & sb & cb \\  \hline
 & ff & 1 & 1 & 0 & 0 & 0 & 0 & 0 \\ \cline{2-9}
 $\pi_F$ &  lb & 0 & 0 & 1 & 0 & 0 & 0 & 1 \\ \cline{2-9}
 &  rb & 0 & 0 & 0 & 1 & 0 & 1 & 0 \\ \hline  
 & ff & 1 & 0 & 1 & 0 & 0 & 0 & 0 \\ \cline{2-9}
$\pi_S$ &  lb & 0 & 0 & 0 & 1 & 1 & 0 & 0 \\ \cline{2-9}
 &  rb & 0 & 1 & 0 & 0 & 0 & 0 & 1 \\ \hline
 & ff & 1 & 0 & 0 & 1 & 0 & 0 & 0 \\ \cline{2-9}
$\pi_C$ &  lb & 0 & 0 & 1 & 0 & 1 & 0 & 0 \\ \cline{2-9}
 &  rb & 0 & 1 & 0 & 0 & 0 & 1 & 0 \\ \hline
 &  $\nu$ & 0 & 0 & 0 & n+2 & 0 & 0 & 0 \\ \hline
\end{tabular}
\end{center}
\vspace*{.2in}

Also in the Table there is a `density row', labelled `$\nu$' which is important
in the description of the composition results. These exponents are fixed by
the natural identification of density bundles:
\[ (\pi_{F})^*K\! D\otimes(\pi_{S})^*K\! D\otimes(\pi_{C})^*(K\! D')
\otimes|dt|^{\hf}\cong\prod\limits_{F=f\! F, f\! f, s\! f, c\! f, fb, sb, cb}
\rho^{\nu_F}_F\cdot \Omega(W_0). \]
Here $K\! D'$ is the half-density bundle with the opposite weighting to $K\! D$
so that $K\! D'\otimes K\! D \cong \Omega.$

 Let $A \in \Psi_{h0}^{-\infty,{\cal E}} (M; \Omega^{\hf}), \ B \in 
\Psi_{h0}^{-\infty,{\cal F}} (M, \Omega^{\hf})$. The composition $C = A \circ B$ 
can be written in terms of their Schwartz kernels via the b-fibrations 
$\pi_{O}, \ O=F, S, C$:
\[ \kappa_C KD' = (\pi_{C})_* [ (\pi_{S})^* \kappa_A \cdot 
(\pi_{F})^* \kappa_B \cdot (\pi_{C})^* (KD') \cdot 
(\tilde{\pi}_{C})^* (|dt'|^{\hf})], \]
where $\tilde{\pi}_{C} = \pi_C^2 \circ \beta_0^3: W_0 \rightarrow Z$ and $\beta_0^3$ is the blow down map.

By the push-forward theorem
\[ C \in \Psi_{h0}^{-\infty,{\cal G}} (M, \Omega^{\hf}) \]
for some index set ${\cal G}$. The index set ${\cal G} = (G_{f\! f}, G_{lb}, 
G_{rb})$ can be computed by the Mellin transform
\[ \langle \kappa_C \cdot KD', \rho_{f\! f}^{z_1} \rho_{lb}^{z_2} 
\rho_{rb}^{z_3}\rangle = 
\langle (\pi_{F})^* (\tilde{\kappa}_A \cdot (\pi_{S})^* 
(\tilde{\kappa}_B) \cdot \prod_{F} \rho_F^{\nu_F} \cdot \Omega(W_0), 
(\pi_{C})^* (\rho_{f\! f}^{z_1} \rho_{lb}^{z_2} \rho_{rb}^{z_3}) \rangle,
\]
where $\kappa_A = \tilde{\kappa}_A \cdot K\! D, \ \kappa_B = \tilde{\kappa}_B
\cdot K\! D$. From the table we obtain the formulas for our index set.
\qed

\subsection{Uniform structure of heat kernel}

The proceeding construction enables us to prove the following theorem, which 
gives the uniform structure of the heat kernel for absolute (or relative) 
boundary condition.
\begin{theo} \label{ushk}
If we denote by $e^{-t\Delta}$ the heat kernel of the Hodge Laplacian for the 
absolute (or relative) boundary condition, then
\[ e^{-t\Delta} \in \Psi_{h0}^{2, 2, 0, 0}(M, \Lambda M\otimes \xi). \]
Moreover the normal operators of $e^{-t\Delta}$ at the front face and the 
temporal face are explicitly computable.
\end{theo}

\Pf Our proof is actually constructive. Namely we construct a (unique) $H\in 
\Psi_{h0}^{I}(M; \Lambda M\otimes \xi)$ where $I=\{ 2, 2, 0, 0\}$
such that 
\be 
t(\pt_{t}+\Delta)H=0 \ \mbox{mod}\ \Psi_{h0}^{I'}(M; \Lambda M\otimes \xi)
\ee
for $I'=\{3, 3, 1, 1\}$, and
\be \label{ic}
N_{h, 2}(H)=\mbox{Id},
\ee
and H satisfies the boundary condition
\[ P_a(H|_{B_{lb}})=0 , \ P_a(dH|_{B_{lb}})=0. \]

Equations (\ref{ushk}) and (\ref{ic}) translate to conditions on the 
normal operators of $H$ by using Proposition~\ref{comp1}:                                                       
\be \label{ctf}
( ^{0}\sigma_2(t\Delta) -\frac{1}{2}(R+n))N_{h,2}(H)=0, \
\int_{\mbox{fiber}} N_{h,2}(H)=\mbox{Id},
\ee
\be \label{cff}
(\pt_{T}+\Delta_E)N_{f,2}(H)=0.
\ee
Finally the boundary condition translates into the same type boundary 
condition for (\ref{cff}).

The first equation is a fiber by fiber differential equation and can be solved 
uniquely subject to the integral condition. Furthermore because of the 
compatibility condition this fixes the integral conditions for (\ref{cff}). 
Thus the normal operator $N_{h, 2}(H)$ is 
necessarily the heat kernel for the elliptic boundary condition on the half 
tangent space.
These two operators have the same indicial family, so using the existence part 
of the compatibility it follows that there is an element $G_1 \in 
\Psi_{h0}^{I}(M; \Lambda M\otimes \xi)$ satisfying the symbolic conditions (\ref{ctf}), 
(\ref{cff}).

This first approximation therefore satisfies 
\be \label{fa}
t(\pt_{t}+\Delta)G_1=- R_1,\  R_1 \in 
\Psi_{h0}^{3,3,1,1}(M; \Lambda M\otimes \xi).
\ee
Thus $G_1$ is already a parametrix. We now modify $G_1$. Using the heat calculus 
we can find a $G_0 \in \Psi_{h0}^{2,2,\infty, \infty}(M, \Lambda M\otimes \xi)$ such that 
\[ t (\pt_t+ \Delta) G_0 = t \Id -t R_0, \]
where $R_0 \in \Psi_{h0}^{0,\infty, \infty, \infty}(M, \Lambda M\otimes \xi)$.

It follows that there is a correction term $G_0' \in \Psi_{h0}$ such that the 
modification of the parametrix $G_2 =G_1 -G_0'$ still has the correct normal 
operator and is a parametrix in the strong sense that
\[ t(\pt_t+ \Delta) G_2 = t \Id -t R_2, \]
where $R_2 \in \Psi_{h0}^{\infty,3, 1,1}$.

By Proposition \ref{comp2}, $(R_2)^k \in  \Psi_{h0}^{\infty,3k,1, 1}$.
Thus the Neumann series $\sum_{k=0}^\infty (R_2)^k$ can be summed modulo a term 
vanishing rapidly at $B_{f\! f}$, i.e. there exists $S' \in \Psi_{h0}^{\infty, 
1,1}$ such that 
\[ (Id - S') (\Id -R_2) = \Id -R_3, \ \ \ R_3 \in \Psi_{h0}^{\infty,\infty,1, 
1}.\]
Thus Id $-R_3$ can be inverted with an operator of the same type. It follows 
that Id $-R_2$ has a two-side inverse Id $-S, \ S \in \Psi_{h0}^{-\infty,1,1} 
(M, F)$. This in turn means we have  $\exp (-t \Delta) = G_1 (Id -S) = G_1 - G_1 
\circ S$
and we obtain our result.
\qed

\medskip
In the discussion below we will be concerned with the trace of the elements of 
$\Psi_{h0}^{k,k}(M; E)
= \Psi_{h0}^{k,k,0,0}(M; E)$. The elements of $\Psi_{h0}^{k,k}(M; E)$ are families of smoothing 
operators on $M$, 
hence trace class. By Lidsky's theorem the trace is the integral over the 
diagonal of the pointwise trace of the kernel. Consider the diagonal $\Delta(M) \times {\Bbb 
R}_+$ of $Z$, which lifts to an embedded submanifold of $M_{h0}^2$. Let $i:\ 
\Delta_{h0} \ra M_{h0}^2$ be the 
embedding of the lifted diagonal. We have
\be \label{ld}
\ \Delta_{h0} \cong M_{h0}=[M\times {\Bbb R}_+; \dM\times\{0\}, S]
\ee
which blows down to $ M\times {\Bbb R}_+$. Restricted to the lifted diagonal, the kernels of  the elements of 
$\Psi_{h0}^{k,k}(M; E)$ can be interpreted as a density:
\[ Hom(E)\otimes\Omega^{1/2}(M_{h0}^2)|_{\Delta_{h0}} \cong hom(E)\otimes 
\Omega(M_{h0})\otimes \Omega^{\hf}({\Bbb R}_+). \]
Thus the trace of $A\in \Psi_{h0}^{2,2,0,0}(M; E)$ is the 
push-forward to ${\Bbb R}_+$ of the density
\[ (\tr\, A)|_{\Delta_{h0}} \in C^{\infty}(M_{h0}; \Omega(M_{h0})\otimes \Omega^{\hf}({\Bbb R}_+)). \]

\begin{lem} \label{mmlem}
As a map 
\[ \Tr:\ \Psi_{h0}^{k,k}(M; E) \ra t^{-\o{n-k+2}{2}} (C^{\infty}({\Bbb R}_+, \Omega^{\hf}({\Bbb R}_+))+ t^{\hf}
C^{\infty}({\Bbb R}_+, \Omega^{\hf}({\Bbb R}_+)) +O(t\log t)).  \]
i.e.
\[ \Tr(A)=t^{-\o{n-k+2}{2}}(r_A(t)+ t^{\hf} s_A(t) + O(t\log t))\,|dt|^{\hf}, \]
where $r_A(t)$, $s_A(t)$ are smooth functions of $t$. Moreover the 
leading terms can be computed by
\be \label{lftr}
 r_A(0)  = \int_M \tr(N_{h,k}(A)), \ \ \  s_A(0) = \int_{B_f} \tr( N_{f,k}(A)),
\ee
where $M \rightarrow\ ^0\!TM$ is the $0$-section and $B_f$ is the front face of $M_{h0}$.
\end{lem}
\Pf This follows from the pushforward theorem, using Mellin transform. Namely, if we denote by $\beta'$ the composition of the blowdown map $\beta_1:\, M_{h0}=[M\times {\Bbb R}_+; \dM \times \{0\}, S] \ra M\times {\Bbb R}_+$ and the projection onto the second factor
$M\times {\Bbb R}_+ \ra {\Bbb R}_+$, we have 
\[ \langle (\beta')_* ( (\tr\, A)|_{\Delta_{h0}}), \ t^z|dt|^{\hf} \rangle =  \langle (\tr\, A)|_{\Delta_{h0}}, \  (\beta')^*(t^z|dt|^{\hf}) \rangle. \]
Now 
\ban
 K\! D |_{\Delta_{h0}} \otimes (\beta')^*(|dt|^{\hf}) & = & (\beta_1)^*(t^{-\o{n+2}{2}} \Omega^{\hf}(M^2\times {\Bbb R}_+) \otimes \Omega^{\hf}({\Bbb R}_+) |_{\Delta(M \times {\Bbb R}_+)}) \\
 & = & \rho_{f\! f}^{-n} \rho_{t\! f}^{-n-1} \Omega (M_{h0}),
\ean
which gives us the desired result.
\qed

\section{Clifford Structure and Cancellation}

In order to prove our main result we need to modify our previous construction by 
incorporating the local index theorem type computation of Bismut-Zhang. A crucial 
ingredient is the Clifford structure of the exterior algebra, which we recall 
now.

\subsection{Clifford Structure of the Exterior Algebra}

If $E$ is a finite dimenional vector space of dimension $n$, the exterior 
algebra $\Lambda(E^*)$ is naturally ${\Bbb Z}$-graded, which induces a natural  
${\Bbb Z}_2$-grading. If $A\in$ End$(\Lambda(E^*))$, we let $\Trs (A)$ denote the supertrace of $A$.

Now assume that $E$ is equipped with an inner product $g$. If $e\in E$, let 
$e^*\in E^*$ correspond to $e$ by the metric $g$. Set
\ban \label{csea}
c(e) & = & e^*\wedge - i_e , \\
\hc(e) & = & e^*\wedge + i_e.
\ean
These define elements of $\End(\Lambda(E^*))$. Moreover, we have the following 
commutator relations
\begin{eqnarray} \label{cr}
c(e)c(e')+c(e')c(e) & = & -2 \langle e, e' \rangle, \nonumber \\
\hc(e)\hc(e')+\hc(e')\hc(e) & = & 2 \langle e, e' \rangle, \\
c(e)\hc(e')+\hc(e')c(e) & = & 0. \nonumber
\end{eqnarray}

{}From (\ref{cr}), we see that $c$ and $\hc$ extend to representations of the 
Clifford algebra of $E$. Furthermore, one verifies that $\End(\Lambda(E^*))$ is 
generated as an algebra by $1$ and the $c(e), \hc(e)$'s.

Let $e_1, \cdots, e_n$ be an orthonormal basis of $E$, and $e^1, \cdots, e^n$ 
the dual basis of $E^*$. A simple but essential algebraic fact is the following 
result which we quote from \cite[Proposition 4.9]{bz}.

\begin{lem} \label{bzlem}
Among the monomials in the $c(e), \hc(e)$'s, only $c(e_1\hsp)\hc(e_1\hsp)
\hsp\cdots c(e_n\hsp)\hc(e_n\hsp)$ has a nozero supertrace, which is 
\[ \Trs [c(e_1)\hc(e_1) \cdots c(e_n)\hc(e_n)]=(-2)^n. \]
\end{lem}

\subsection{Supertrace and Cancellation}

The above discussion applied fiberwise to $TM$ shows that $\End(\Lambda^*M)$
is generated as an algebra by $1$ and the $c(e), \hc(e)$'s  for any local 
orthonormal basis $e$. This gives a natural filtration of $\End(\Lambda^*M)$, 
which we will exploit by making a global rescaling of $\End(\Lambda^*M\otimes 
\xi)$ near the front and temporal faces of $M^2_{h0}$. Since this is localized 
near the diagonal, the vector bundle $\xi$ does not appear in the discussion and we happily 
suppress its presence. 

Denote by
\[ F^0 \subset F^1 \subset \cdots \subset F^n \]
the natural filtration on $\End(\Lambda^*M)$ induced by the Clifford structure introduced above, i.e., $F^j$ is the subspace generated by monomials of length 
$\leq j$ in $c(e), \hc(e)$.  Its associated graded algebra $G=\oplus G^i$, $G^i = F^i/F^{i-1}$, is naturally 
isomorphic to $\Lambda(T^*M) \otimes \Lambda(T^*M)$ (Cf. (\ref{cr})). We distinguish the elements of the second copy
of $\Lambda(T^*M)$ from the first by putting a hat on them. For example ${\cal R}=\o{1}{8} R_{ijkl}\, e^i 
\wedge e^j \wedge \hat{e}^k \wedge \hat{e}^l \in \Lambda(T^*M) \otimes \Lambda(T^*M)$.

Recall that the Berezin integral \cite{bz} takes elements of $\Lambda(T^*M) \otimes \Lambda(T^*M)$ and sends them 
to $\Lambda(T^*M)$:
\be \int^B :\  \Lambda(T^*M) \otimes \Lambda(T^*M) \longrightarrow \Lambda(T^*M)  \label{bint}
\ee
The Berezin integral is intimately related to the supertrace and cancellation mechanism.

Applying the construction in \cite{dm} we obtain a rescaled version of the 
pseudodifferential \ calculus \ constructed \ above, 
\ $\Psi_{h0, G}^I(M; \Lambda^*M\otimes \xi)$. \ For 
\ example, $\Psi_{h0, G}^{2,2,0,0}(M; \Lambda^*M\otimes 
\xi)$ consists of those elements of 
$\Psi_{h0}^{2,2,0,0}(M; \Lambda^*M\otimes 
\xi)$ which have the following asymptotic expansion near 
either the front or temporal face $F$:
\be \label{resc}
k=\rho_{f\! f}^2 \rho_{t\! f}^2 \sum_{j=0}^{n}\rho_F^{\o{j}{2}} k_j, \ \ \ k_j \in C^{\infty}(M^2_{h0}, 
F^j\otimes K\! D).
\ee
The normal operators $N_{h,k, G}, \ N_{f,k, G}$ of elements of  $\Psi_{h0, G}^I(M; \Lambda^*M\otimes 
\xi)$ will now be taking values in the associated graded algebra $G \cong  \Lambda(T^*M) \otimes \Lambda(T^*M)$.

The advantage of the modification comes from the following cancellation result.
(Note the disappearance of the singular factor $ t^{-\o{n}{2}}$.)
\begin{lem} \label{mlem}
As a map 
\[ \Trs:\ \Psi_{h0, G}^{k,k,0,0}(M; \Lambda^*M\otimes 
\xi) \ra  t^{\o{k-2}{2}} (C^{\infty}({\Bbb R}_+, \Omega^{\hf}({\Bbb R}_+))
+t^{\hf} C^{\infty}({\Bbb R}_+, \Omega^{\hf}({\Bbb R}_+)) + O(t \log t)).  \]
i.e.
\[ \Trs(A)= t^{\o{k-2}{2}} (r_A(t) +t^{\hf} s_A(t)  + O(t \log t))\, |dt|^{\hf}, \]
where $r_A(t)$, $s_A(t)$ are smooth functions of $t$. Moreover the 
leading terms can be computed by
\be \label{glftr}
 r_A(0)  =  \int_M \int^B N_{h,k, G}(A), \ \ \ s_A(0) = \int_{B_f} \int^B N_{f,k, G}(A). 
\ee
\end{lem}

\Pf This comes from the definition and the main algebraic fact Lemma \ref{bzlem}.
\qed

\medskip
With this in mind we have the following significant refinement on the structure 
of the heat kernel.

\begin{theo} \label{gushk}
 The heat kernel $e^{-t\Delta}$ of the Hodge Laplacian for the absolute (or 
relative) boundary condition, is an element of the rescaled heat calculus:
\[ e^{-t\Delta} \in \Psi_{h0, G}^{2, 2, 0, 0}(M, \Lambda^*M\otimes 
\xi). \]
Moreover the normal operators of $e^{-t\Delta}$ at the front face and the 
temporal face are given by
\ban \label{gno}
 N_{h,2, G}& = & (4\pi)^{-\o{n}{2}}\exp(-|v|^2)\exp(-{\cal R}), \label{gnoah}\\
N_{f,2, G} & = & \exp(-T(\Delta_E+{\cal R})), \label{gnoaf}
\ean
where ${\cal R}=\o{1}{8} R_{ijkl}\, e^i 
\wedge e^j \wedge \hat{e}^k \wedge \hat{e}^l$ and $\Delta_E$ is the fibrewise 
Laplacian on the inward pointing half of tangent bundle $TM$ with the absolute 
(relative respectively) boundary condition.
\end{theo}

\Pf We only need to adapt the previous argument to the rescaling. 
According to the Bochner-Weitzenb\"ock formula,
\[ \Lap  =  \cLap-R_{ijkl}\ext(e_i)\inn(e_j) \ext(e_k)\inn(e_l). \]
Here $\cLap=\nabla^*\nabla$ is the connection Laplacian; with respect to any local
orthonormal frame of $TM$ it is
\[ \cLap=-\sum_{i}(\nabla^2_{e_{i}}-\nabla_{\nabla_{e_i}e_i}). \]
Expressing everything in terms of $c(e)$ and $\hc(e)$, we have
\[ \Lap  =  \cLap +  c(R) -\o{1}{4} S, \]
with $S$ the scalar curvature and 
\[ c(R)=\o{1}{8}  R_{ijkl}\, c(e_i) c(e_j) \hc(e_k) \hc(e_l). \]

Now since each Clifford element carries a weight of $\rho^{\hf}$ at both the temporal and the front faces, it follows then for $A\in \Psi_{h0, G}^I(M; \Lambda^*M\otimes \xi)$,  $tc(R)\circ A \in \Psi_{h0, G}^I(M; \Lambda^*M\otimes 
\xi)$, and 
\[  N_{h, G}(tc(R)\circ A)= {\cal R}N_{h, G}(A), \ \ \  N_{f, G}(tc(R)\circ A)= T{\cal R}N_{f, G}(A). \]
Finally, for the connection, we note that
\[ \nabla_{e_i} = \Gamma_{ij}^k \ext(e_k) \inn(e_j) = \o{1}{4} \Gamma_{ij}^k (c(e_k)+\hc(e_k))(\hc(e_j)-c(e_j)) \]
again produces only lower order terms.
\qed

\subsection{The Model Problem}

Our model problem is to solve in ${\Bbb R}^n_+$ for the heat kernel of $\Delta 
+{\cal R}$, where 
${\cal R}=\o{1}{8} R_{ijkl}\, e^i 
\wedge e^j \wedge \hat{e}^k \wedge \hat{e}^l$, 
subject to the, say, absolute boundary conditions. More specifically, let $(x,\, 
y) \in {\Bbb R}_+ \times {\Bbb R}^{n-1} = {\Bbb R}^n_+$ and $e_0=\pt_x$. Then we
need to solve

\be 
\left\{ \begin{array}{l} (\pt_t +\Delta +{\cal R})K=0  \ \mbox{in}\ {\Bbb R}^n_+ 
\\
K|_{t=0}=\mbox{I} \\
e^0\wedge i_{e_0} K|_{x=0}=0 \\
i_{e_0} e^0 \wedge \o{\pt K}{\pt x}|_{x=0}=0.  \end{array} \right.
\ee

\begin{lem} \label{comm} Let ${\cal R}_0= \o{1}{4}R_{0jkl}\, e^0 \wedge e^j 
\wedge \hat{e}^k \wedge \hat{e}^l$. We have
\be \left[ {\cal R}, \ i_{e_0} e^0 \wedge \right] = {\cal R}_0, \ \ 
\left[ {\cal R}, \ e^0 \wedge i_{e_0} \right] = -{\cal R}_0.
\ee
\end{lem}

\Pf This is a straightforward computation using
\[ i_{e_0} e^0 \wedge +e^0 \wedge i_{e_0}=1.  \]
In fact, 
\begin{eqnarray*} {\cal R} (i_{e_0} e^0 \wedge) & = & 
\left( \o{1}{8}\sum_{i, j\neq 0} R_{ijkl}\, e^i \wedge e^j \wedge \hat{e}^k 
\wedge \hat{e}^l + {\cal R}_0 \right) (i_{e_0} e^0 \wedge) \\
& = & (i_{e_0} e^0 \wedge) {\cal R} +  {\cal R}_0 
\end{eqnarray*}

The other formula is obtained similarly. \qed

\medskip
Now, let 
\begin{eqnarray*} K_D(t) & = & \o{1}{(4\pi t)^{n/2}} (e^{-(x-x')^2/4t} - 
e^{-(x+x')^2/4t}) e^{-|y-y'|^2/4t}, \\
 K_N(t) & = & \o{1}{(4\pi t)^{n/2}} (e^{-(x-x')^2/4t} + e^{-(x+x')^2/4t}) 
e^{-|y-y'|^2/4t}
 \end{eqnarray*} 
be the Dirichlet and Neumann heat kernels respectively. Also, for $k=k(t,x, 
y, x', y')$, $\tilde{k}=\tilde{k}(t,\, x, \, y, \, x', \, y')$, their 
convolution is
\be k * \tilde{k} = \int_0^t \int_0^{\infty} \int_{{\Bbb R}^{n-1}} k(s,\, x, \, 
y, \, x', \, y') \tilde{k}(t-s,\, x', \, y', \, x'', \, y'') \, dt \, dx' \, 
dy'.
\ee

\begin{prop} Let \be K_0= i_{e_0} e^0 \wedge K_N(t) e^{-t{\cal R}} + e^0 \wedge
i_{e_0} K_D(t) e^{-t{\cal R}}.  \ee Then \be (\pt_t +\Delta +{\cal R})K_0= {\cal
R}_0 [K_N(t)-K_D(t)]e^{-t{\cal R}}\equiv K_1.  \ee And therefore, the solution 
of
the model problem is \be K=K_0 + K_0 * K_1.  \ee \end{prop}

\Pf The first equation follows immediately from Lemma~\ref{comm}.  Then we
note that $K_0$ satisfies the absolute boundary condition and we apply the
Duhamel's principle.  The infinite sum terminates at the second term since $
{\cal R}_0$ contains the nilpotent element $e^0 \wedge$.  \qed

\subsection{Proof of the Theorem}

Finally we are in a position to prove Theorem~\ref{mt}. 

{\em Proof of Theorem~\ref{mt}}:  As we mentioned before, we use the metric deformation (\ref{md}) in the variation formula (\ref{vf}). To compute the right hand side of (\ref{vf}) we  
rewrite it as
\be \label{rweta}
\lim_{t \rightarrow 0} \Trs F, \ 
F=*_l^{-1}\o{\pt *_l}{\pt l}e^{-t\Delta_l}.
\ee
Since
\begin{eqnarray}
*_l^{-1}\o{\pt *_l}{\pt l} &  =  & - \sum_{1\leq i, j \leq n}\hf \langle (g_l^{-1} \o{\pt g_l}{\pt l})(e_i), \, e_j \rangle c(e_i) \hc(e_j),  \label{star} \\
& = & - h_{ab} c(e_a) \hc(e_b) x + O(x^2).
\end{eqnarray}
(see \cite[Proposition 4.15]{bz}),  from 
Theorem~\ref{gushk} (and Theorem~\ref{ushk}), we have
\be \label{poly}
F\in \Psi_{h0, G}^{1,1,0,0}(M; E).
\ee
Therefore, according to Lemma~\ref{mlem},
\[ \Trs(F)= t^{-\o{1}{2}} (r(t) +t^{\hf} s(t))\, |dt|^{\hf}, \]
with $r(0)$ and $s(0)$ given by (\ref{glftr}). 

Thus we must compute the two terms 
in (\ref{glftr}). The first term computes the contribution from the interior of 
the manifolds and, just as in the case of the closed manifolds, can be seen to 
be zero. 
The term that really contributes to the variation formula is the second term which computes the boundary contribution. That is where 
the model problem comes in.

To compute the contribution of the model problem, we need to evaluate $K$ on the spatial  diagonal and at $t=1$, and compute the corresponding Berezin integral, then 
integrate finally over ${\Bbb R}_+ \times \pt M$. Restricted to the spatial 
diagonal and $t=1$, we have
\[ K= \o{1}{(4\pi)^{n/2}} e^{-{\cal R}} + \o{1}{(4\pi)^{n/2}} e^{-x^2} c(e_0) \hc(e_0) 
e^{-{\cal R}} + \o{1}{(4\pi)^{(n-1)/2}} f(x) {\cal R}_0 e^{-{\cal R}},\]
where 
\be f(x)=- \int_0^1 \int_0^{\infty} (e^{-(x-x')^2/4s} + e^{-(x+x')^2/4s})
2 e^{-(x+x')^2/4(1-s)}\, dx'\, ds.
\ee
Thus, we obtain
\ban
s(0) & = & -\,{\rm rank}\,\rho \int_{{\Bbb R}_+ \times \pt M} \int^B h_{ab} e^a \wedge \he^b x \o{1}{(4\pi)^{n/2}} e^{-x^2} e^0 \wedge \he^0 \wedge 
e^{-{\cal R}}  \\
& & - {\rm rank}\,\rho \int_{{\Bbb R}_+ \times \pt M} \int^B h_{ab} e^a \wedge \he^b x \o{1}{(4\pi)^{(n-1)/2}} f(x) {\cal R}_0 e^{-{\cal R}}. 
\ean
The first Berezin integral here can be explicitly evaluated: 
\begin{eqnarray*}
& & -\int_{{\Bbb R}_+ \times \pt M} \int^B h_{ab} e^a \wedge \he^b x \o{1}{(4\pi)^{n/2}} e^{-x^2}  
e^0 \wedge \he^0 \wedge e^{-{\cal R}} \\ 
 &  = &\int_{{\Bbb R}_+ \times \pt M} \o{\pt}{\pt l} i^*(\te(g_0, g_l))  e^{-x^2} xdx  \\
 &  = & \o{\pt}{\pt l}\int_{\pt M} i^*(\te(g_0, g_l)). 
\end{eqnarray*}
Taking $t$ to $0$ and then integrating $l$ from $0$ to $1$ finally gives us the
desired result.
\qed

\def\em{\it}

\end{document}